\documentclass[a4paper,oneside]{amsart}

\RequirePackage{amsmath}
\RequirePackage{bm}
\RequirePackage{amssymb}
\RequirePackage{upref}
\RequirePackage{amsthm}
\RequirePackage{enumerate}
\RequirePackage{pb-diagram}
\RequirePackage{amsfonts}
\RequirePackage[mathscr]{eucal}
\RequirePackage{verbatim}
\RequirePackage{xr}
\RequirePackage{graphicx}
\RequirePackage{floatflt}
\usepackage{calc}
\usepackage{xspace}

\newcommand{\cf}{cf.\@\xspace}
\newcommand{\resp}{resp.\@\xspace}

\newcommand{\al}{\alpha}
\newcommand{\bet}{\beta}
\newcommand{\ga}{\gamma}
\newcommand{\de}{\delta }
\newcommand{\e}{\epsilon}

\newcommand{\f}{\varphi}
\newcommand{\h}{\eta}

\newcommand{\ka}{\kappa}
\newcommand{\lam}{\lambda}

\newcommand{\n}{\nu}

\newcommand{\s}{\sigma}
\newcommand{\x}{\xi}

\newcommand{\C}{\varGamma}
\newcommand{\D}{\varDelta}
\newcommand{\F}{\varPhi}
\newcommand{\Lam}{\varLambda}
\newcommand{\Om}{\varOmega}

\newcommand{\di}[1]{#1\nobreakdash-\hspace{0pt}dimensional}
\newcommand{\nbdd}{\nobreakdash--}
\newcommand{\nbd}{\nobreakdash-\hspace{0pt}}

\newcommand{\fmo}[1]{F_{|_{#1}}}

\newcommand{\fv}[2]{#1\hspace{0pt}_{|_{#2}}}

\newcommand{\const}{\tup{const}}
\newcommand{\ndash}{\nobreakdash--}

\newcommand{\pih}{\frac{\pi}{2}}

\newcommand{\msp[1]}[1]{\mspace{#1mu}}

\newcommand{\R}[1][n+1]{{\protect\mathbb R}^{#1}}

\newcommand{\N}{{\protect\mathbb N}}

\newcommand{\eR}{\stackrel{\lower1ex \hbox{\rule{6.5pt}{0.5pt}}}{\msp[3]\R[]}}
\newcommand{\eN}{\stackrel{\lower1ex \hbox{\rule{6.5pt}{0.5pt}}}{\msp[1]\N}}
\newcommand{\eO}{\stackrel{\lower1ex
\hbox{\rule{6pt}{0.5pt}}}{\msc O}}

\DeclareMathOperator{\graph}{graph}

\DeclareMathOperator{\id}{id}

\newcommand\im{\implies}
\newcommand\ra{\rightarrow}

\newcommand\hra{\hookrightarrow}

\newcommand\pa{\partial}
\newcommand\pde[2]{\frac {\partial#1}{\partial#2}}
\newcommand\pd[3]{\frac {\partial#1}{\partial#2^#3}}   
 
\newcommand\df[2]{\frac {d#1}{d#2}}

\newcommand{\un}{\infty}
\newcommand{\A}{\forall}

\newcommand{\set}[2]{\{\,#1\colon #2\,\}}
\newcommand{\uu}{\cup}
\newcommand{\ii}{\cap}
\newcommand{\uuu}{\bigcup}

\newcommand{\uud}{ \stackrel{\lower 1ex \hbox {.}}{\uu}}
\newcommand{\uuud}[1]{ \stackrel{\lower 1ex \hbox {.}}{\uuu_{#1}}}
\newcommand\su{\subset}
\newcommand\Su{\Subset}

\newcommand\eS{\emptyset}
\newcommand{\sminus}[1][28]{\raise 0.#1ex\hbox{$\scriptstyle\setminus$}}
\newcommand{\cpl}{\complement}

\newcommand{\wed}{\wedge}

\newcommand\ti{\times }

\newcommand{\abs}[1]{\lvert#1\rvert}

\newcommand{\norm}[1]{\lVert#1\rVert}

\newcommand{\spd}[2]{\protect\langle #1,#2\protect\rangle}

\newcommand\ch[3]{\varGamma_{#1#2}^#3}
\newcommand\cha[3]{{\bar\varGamma}_{#1#2}^#3}

\newcommand{\riem}[4]{R_{#1#2#3#4}}
\newcommand{\riema}[4]{{\bar R}_{#1#2#3#4}}

\newcommand{\tit}{\textit}

\newcommand{\tup}{\textup}

\newcommand{\mc}{\protect\mathcal}
\newcommand{\msc}{\protect\mathscr}

\providecommand{\bysame}{\makebox[3em]{\hrulefill}\thinspace}

\newcommand{\ci}{\cite}

\newcommand{\cq}[1]{\glqq{#1}\grqq\,}

\newcommand{\bt}{\begin{thm}}
\newcommand{\bl}{\begin{lem}}
\newcommand{\bc}{\begin{cor}}
\newcommand{\bd}{\begin{definition}}
\newcommand{\bpp}{\begin{prop}}
\newcommand{\br}{\begin{rem}}
\newcommand{\bn}{\begin{note}}
\newcommand{\be}{\begin{ex}}
\newcommand{\bes}{\begin{exs}}
\newcommand{\bb}{\begin{example}}
\newcommand{\bbs}{\begin{examples}}
\newcommand{\ba}{\begin{axiom}}
\newcommand{\bas}{\begin{assumption}}

\newcommand{\et}{\end{thm}}
\newcommand{\el}{\end{lem}}
\newcommand{\ec}{\end{cor}}
\newcommand{\ed}{\end{definition}}
\newcommand{\epp}{\end{prop}}
\newcommand{\er}{\end{rem}}
\newcommand{\en}{\end{note}}
\newcommand{\ee}{\end{ex}}
\newcommand{\ees}{\end{exs}}
\newcommand{\eb}{\end{example}}
\newcommand{\ebs}{\end{examples}}
\newcommand{\ea}{\end{axiom}}
\newcommand{\eas}{\end{assumption}}

\newcommand{\bp}{\begin{proof}}
\newcommand{\ep}{\end{proof}}
\newcommand{\eps}{\renewcommand{\qed}{}\end{proof}}

\newcommand{\bal}{\begin{align}}

\newcommand{\bi}[1][1.]{\begin{enumerate}[\upshape #1]}
\newcommand{\bia}[1][(1)]{\begin{enumerate}[\upshape #1]}
\newcommand{\bin}[1][1]{\begin{enumerate}[\upshape\bfseries #1]}
\newcommand{\bir}[1][(i)]{\begin{enumerate}[\upshape #1]}
\newcommand{\bic}[1][(i)]{\begin{enumerate}[\upshape\hspace{2\cma}#1]}
\newcommand{\bis}[2][1.]{\begin{enumerate}[\upshape\hspace{#2\parindent}#1]}
\newcommand{\ei}{\end{enumerate}}

\newcommand\ndots{\raise 0.47ex \hbox {,}\hskip0.06em\cdots %
     \raise 0.47ex \hbox {,}\hskip0.06em} 

\newcommand{\q}{\quad}
\newcommand{\qq}{\qquad}

\newcommand{\hp}{\hphantom}

\newcommand\nd{\noindent}

\newskip\Csmallskipamount                                                
\Csmallskipamount=\smallskipamount
\newskip\Cmedskipamount
\Cmedskipamount=\medskipamount
\newskip\Cbigskipamount
\Cbigskipamount=\bigskipamount

\newcommand\cvs{\vspace\Csmallskipamount}   
\newcommand\cvm{\vspace\Cmedskipamount}

\newskip\csa
\csa=\smallskipamount

\newskip\cma
\cma=\medskipamount

\newskip\cba
\cba=\bigskipamount

\newdimen\spt
\newcommand\citem{\cvs\advance\itemno by
1{(\romannumeral\the\itemno})\hskip3pt}
\newcommand{\bitem}{\cvm\nd\advance\itemno by
1{\bf\the\itemno}\hspace{\cma}}

\newcount\itemno
\newcommand{\las}[1]{\label{S:#1}}

\newcommand{\lae}[1]{\label{E:#1}}
\newcommand{\lat}[1]{\label{T:#1}}
\newcommand{\lal}[1]{\label{L:#1}}

\newcommand{\lar}[1]{\label{R:#1}}

\newcommand{\laas}[1]{\label{Ass:#1}}

\newcommand{\rs}[1]{Section~\ref{S:#1}}

\newcommand{\rt}[1]{Theorem~\ref{T:#1}}
\newcommand{\rl}[1]{Lemma~\ref{L:#1}}

\newcommand{\ras}[1]{Assumption~\ref{Ass:#1}}

\newcommand{\rr}[1]{Remark~\ref{R:#1}}

\newcommand{\re}[1]{\eqref{E:#1}}

\newskip\thmskip
\thmskip=\parindent

\newskip\hsk
\newenvironment{hinw}{\labelsep=0pt\begin{list}{}{\labelsep=0pt\itemindent=0pt\labelwidth=0pt\leftmargin=\parindent\rightmargin=0pt\partopsep=\cba}%
\item\it\nopagebreak\nopagebreak}%
{\end{list}}

\newcommand\bh{\begin{hinw}}
\newcommand{\eh}{\end{hinw}}

\newtheoremstyle{normal}
  {\cba}
  {\cba}
  {}
  {\thmskip}
  {\bfseries}
  {.}
  {\hsk}
  {}

\newtheoremstyle{abschnitt}
  {\cba}
  {\cba}
  {}
  {\thmskip}
  {\bfseries}
  {.}
  {\hsk}
  {}

\newtheoremstyle{italic}
  {\cba}
  {\cba}
  {\itshape}
  {\thmskip}
  {\bfseries}
  {.}
  {\hsk}
  {}

\newtheoremstyle{aufgaben}
  {\cba}
  {\cba}
  {}
  {}
  {\normalsize\bfseries}
  {.}
  {\hsk}
  {}

\newtheoremstyle{break}
  {\cba}
  {\cba}
  {\itshape}
  {}
  {\bfseries}
  {.}
  {\newline}
  {}

\swapnumbers
\theoremstyle{italic}
\newtheorem{thm}[subsection]{Theorem}
\newtheorem{lem}[subsection]{Lemma}
\newtheorem{prop}[subsection]{Proposition}
\newtheorem{cor}[subsection]{Corollary}

\theoremstyle{normal}
\newtheorem{rem}[subsection]{Remark}
\newtheorem{definition}[subsection]{Definition}
\newtheorem{example}[subsection]{Example}
\newtheorem{examples}[subsection]{Examples}
\newtheorem{ex}[subsection]{Exercise}
\newtheorem{note}[subsection]{}
\newtheorem{axiom}[subsection]{Axiom}
\newtheorem{assumption}[subsection]{Assumption}

\theoremstyle{aufgaben}
\newtheorem{exs}[subsection]{Exercises}

\swapnumbers

\numberwithin{equation}{section}
\numberwithin{figure}{section}

\newenvironment{textequation}[1][0.8]
{\begin{equation}
\begin{aligned}
\begin{minipage}{#1\linewidth}}
{\end{minipage}
\end{aligned}
\end{equation}
\ignorespacesafterend}

\newcommand{\btext}{\begin{textequation}}
\newcommand{\etext}{\end{textequation}}

\usepackage[german,english]{babel}
\usepackage{graphicx}
\RequirePackage{amsmath}
\RequirePackage{bm}
\RequirePackage{amssymb}
\RequirePackage{upref}
\RequirePackage{amsthm}
\RequirePackage{enumerate}
\RequirePackage{pb-diagram}
\RequirePackage{amsfonts}
\RequirePackage[mathscr]{eucal}
\RequirePackage{verbatim}
\RequirePackage{xr}
\RequirePackage{graphicx}
\usepackage{calc}
\usepackage{xspace}

\RequirePackage{upref}
\RequirePackage{amsthm}
\RequirePackage{enumerate}
\usepackage[mathscr]{eucal}

\usepackage{xr-hyper}

\usepackage{calc}

\newlength{\oddsidemarginlength}
\newlength{\topmarginlength}

\newcounter{numberoflines}
\newcounter{tempcc}
\oddsidemargin=\oddsidemarginlength
\evensidemargin=\oddsidemargin
\topmargin=\topmarginlength
\begin{document}

\flushbottom


\title[Minkowski type problems]{Minkowski type problems for convex hypersurfaces in the sphere}

\author{Claus Gerhardt}
\address{Ruprecht-Karls-Universit\"at, Institut f\"ur Angewandte Mathematik,
Im Neuenheimer Feld 294, 69120 Heidelberg, Germany}
\email{gerhardt@math.uni-heidelberg.de}
\urladdr{http://www.math.uni-heidelberg.de/studinfo/gerhardt/}
\thanks{}

%
\subjclass[2000]{35J60, 53C21, 53C44, 53C50, 58J05}
\keywords{Minkowski problem, hypersurfaces of prescribed curvature, Weingarten hypersurfaces}
\date{\today}
%

\dedicatory{Dedicated to Leon Simon on the occasion of his 60th birthday}

\begin{abstract}
We consider the problem $F=f(\nu)$ for strictly convex, closed hypersurfaces in $S^{n+1}$ and solve it for curvature functions $F$ the inverses of which are of class $(K)$.
\end{abstract}

\maketitle

\tableofcontents

\setcounter{section}{-1}
\section{Introduction}

In the classical Minkowski problem in $\R$ one wants to find a strictly convex closed hypersurface $M\su \R$ such that its Gau{\ss} curvature $K$ equals a given function $f$ defined in the normal space of $M$ or equivalently defined on $S^n$
\begin{equation}
\fv KM=f(\nu).
\end{equation}

The problem has been partially solved by Minkowski \cite{minkowski:minkowski}, Alexandrov \cite{alexandrov:minkowski}, Lewy \cite{lewy:minkowski}, Nirenberg \cite{nirenberg:minkowski}, and Pogorelov \cite{pogorelov:minkowski}, and in full generality by Cheng and Yau \cite{cheng-yau:minkowski}.

Instead of prescribing the Gaussian curvature other curvature functions $F$ can be considered, i.e., one studies the problem
\begin{equation}\lae{0.2}
\fv FM=f(\nu).
\end{equation}

If $F$ is one of the symmetric polynomials $H_k$, $1\le k\le n$, this problem has recently been solved by Guan and Guan \cite{guan:annals}.  They proved that \re{0.2} has a solution, if $f$ is invariant with respect to a fixed point free group of isometries of $S^n$. 

In this paper we consider the problem \re{0.2} for strictly convex hypersurfaces $M\su S^{n+1}$ and for curvature functions $F$ the inverses of which are of class $(K)$, see \rs{1} or \cite[Definition 1.3]{cg:indiana}. These $F$ include all $H_k$, $1\le k\le n$, $\abs{A}^2$, and also any symmetric, convex curvature function homogeneous of degree $1$, \cf \cite[Lemma 1.6]{cg97}.

We shall show in \rs{2} that for any closed strictly convex hypersurface $M\su S^{n+1}$ there exists a Gau{\ss} map
\begin{equation}
x\in M\ra \tilde x\in M^*,
\end{equation}
where $M^*$ is the polar set of $M$. $M^*$ is also strictly convex, as smooth as $M$, and the Gau{\ss} map is a diffeomorphism.

If we consider $M$ as an embedding in $\R[n+2]$ of codimension $2$, so that the tangent spaces $T_x(M)$ and $T_x(S^{n+1})$ can be identified with  subspaces of $T_x(\R[n+2])$, then the image of  the point $x$ under the Gau{\ss} map is exactly the normal vector $\nu\in T_x(S^{n+1})$
\begin{equation}
\tilde x=\nu\in T_x(S^{n+1})\su T_x(\R[n+2]).
\end{equation}

Thus, the equation \re{0.2} can also be written in the form
\begin{equation}\lae{0.5}
\fv FM=f(\tilde x)\qq\A\,Ê x\in M,
\end{equation}
where $f$ is given as a function defined in $S^{n+1}$.

We shall also prove that \re{0.5} has a dual problem, namely,
\begin{equation}
\fv {\tilde F}{M^*}=f^{-1}(\tilde x)\qq\A\,  \tilde x\in M^*,
\end{equation}
where $\tilde F$ is the inverse of $F$
\begin{equation}
\tilde F(\ka_i)=\frac1{F(\ka_i^{-1})}.
\end{equation}

In the dual problem the curvature is not prescribed by a function defined in the normal space, but by a function defined on the hypersurface.

Both problems are equivalent, solving one also leads to a solution of the dual one; notice also that
\begin{equation}
M^{**}=M\q\wed\q \tilde{\tilde x}=x.
\end{equation}

To find a solution we either impose some symmetry requirement with respect to a group of isometries or we assume the existence of barriers.

\bas\laas{0.1}
(i) Let  $G\su O(n+2)$ be a group of orthogonal transformations with a common fixed point $x_0\in S^{n+1}$ and assume that the induced group of isometries in $S^n$, i.e., the equator of the hemisphere with center in $x_0$, is fixed point free.

\cvm
(ii) Let $0<f\in C^5(S^{n+1})$ be invariant with respect to the group $G$, i.e.,
\begin{equation}
f(Ax)=f(x)\qq\A\, x\in S^{n+1},\;\A\, A\in G.
\end{equation}
\eas

Then we shall prove

\bt
Let $F\in C^5(\C_+)$ be a symmetric, positively homogeneous and monotone curvature function such that its inverse $\tilde F$ is of class $(K)$, then the dual problems
\begin{equation}
\fv FM=f(\tilde x)
\end{equation}
and
\begin{equation}\lae{0.11}
\fv{\tilde F}{M^*}=f^{-1}(\tilde x)
\end{equation}
have strictly convex solutions $M$ \resp $M^*$ of class $C^{4,\al}$, $0<\al<1$, such that the hypersurfaces $M$ \resp $M^*$ are invariant with respect to the group $G$. Furthermore, $-x_0$ is an interior point of the convex body $\hat M$ and $x_0$ an interior point of the convex body $\hat{M}^*$ of $M^*$. The convex bodies $\hat M$, $\hat M^*$ are strictly contained in the corresponding open hemispheres $\mc H(-x_0)$ \resp $\mc H(x_0)$.
\et

Instead of imposing some symmetry assumption, a barrier condition will also work.

\bas\laas{0.3} 
Let $M_i$, $i=1,2$, be strictly convex hypersurfaces of class $C^{6,\al}$ contained in an open hemisphere $\mc H(-x_0)$. $M_1$ is said to be a \tit{lower barrier} for the pair $(F,f)$, if
\begin{equation}
\fv F{M_1}\le f,
\end{equation}
and $M_2$ is called an \tit{upper barrier} for $(F,f)$, if
\begin{equation}
\fv F{M_2}\ge f,
\end{equation}
where in both cases the right-hand side $f$ may either  depend on $x\in M_i$ or $\nu\in T_x(S^{n+1})$ for $x\in M_i$, or, in the latter case, equivalently on $\tilde x\in M_i^*$.
\eas

\bt\lat{0.4}
Let $F\in C^5(\C_+)$ be a symmetric, positively homogeneous and monotone curvature function such that its inverse $\tilde F$ is of class $(K)$, let $0<f\in C^5(S^{n+1})$, and assume that there exist upper and lower barriers for $(F,f)$ in the hemisphere $\mc H(-x_0)$ as defined in the \ras{0.3}, where in addition the barriers $M_i$ should bound a connected open set $\Om$ such that the mean curvature vector of $M_1$ should point to the exterior of $\Om$ and the mean curvature vector of $M_2$ should point into $\Om$. Then the dual problems
\begin{equation}
\fv FM=f(\tilde x)
\end{equation}
and
\begin{equation}
\fv{\tilde F}{M^*}=f^{-1}(\tilde x)
\end{equation}
have strictly convex solutions $M$ \resp $M^*$ of class $C^{6,\al}$, $0<\al<1$, such that the  convex bodies $\hat M$, $\hat M^*$ are strictly contained in the open hemispheres $\mc H(-x_0)$ \resp $\mc H(x_0)$.
\et

 The paper is organized as follows: \rs{1} gives an overview of the definitions and conventions we rely on, while the dual relationship between $M$ and $M^*$ and the properties of the Gau{\ss} map are derived in \rs{2}. The curvature estimates are proved in \rs{3}, the lower order estimates in \rs{4}. The next two sections contain a uniqueness result for invariant convex hypersurfaces with constant $F$, and the existence proof in the invariance case, which is based on a  continuity method using Smale's infinite dimensional version of Sard's theorem  \cite{smale:sard}. Finally, in  \rs 7 we prove \rt{0.4}.
 
 \br\lar{0.5}
 Let us emphasize that after \rs{2} we shall only consider equation \re{0.11}. In order to simplify notation we then shall drop the tilde and the other embellishments and shall solve the equation
 \begin{equation}
\fv FM=f(x)
\end{equation}
for a curvature function $F$ of class $(K)$, where we note that we also replaced $f^{-1}$ by $f$. 
 \er
 
 \section{Notations and definitions}\las 1

The main objective of this section is to state the equations of Gau{\ss}, Codazzi,
and Weingarten for  hypersurfaces $M$ in a \di {(n+1)} Riemannian
manifold
$N$.  Geometric quantities in $N$ will be denoted by
$(\bar g_{ \al \bet})$, $(\riema  \al \bet \ga \de)$, etc., and those in $M$ by $(g_{ij}), 
(\riem ijkl)$, etc. Greek indices range from $0$ to $n$ and Latin from $1$ to $n$; the
summation convention is always used. Generic coordinate systems in $N$ resp.
$M$ will be denoted by $(x^ \al)$ resp. $(\x^i)$. Covariant differentiation will
simply be indicated by indices, only in case of possible ambiguity they will be
preceded by a semicolon, i.e., for a function $u$ in $N$, $(u_ \al)$ will be the
gradient and
$(u_{ \al \bet})$ the Hessian, but e.g., the covariant derivative of the curvature
tensor will be abbreviated by $\riema  \al \bet \ga{ \de;\e}$. We also point out that
\begin{equation}
\riema  \al \bet \ga{ \de;i}=\riema  \al \bet \ga{ \de;\e}x_i^\e
\end{equation}
with obvious generalizations to other quantities.

Let $M$ be a  $C^2$-hypersurface 
with  normal $\n$.

In local coordinates, $(x^ \al)$ and $(\x^i)$, the geometric quantities of the
hypersurface $M$ are connected through the following equations
\begin{equation}\lae{1.2}
x_{ij}^ \al= -h_{ij}\n^ \al
\end{equation}
the so-called \tit{Gau{\ss} formula}. Here, and also in the sequel, a covariant
derivative is always a \tit{full} tensor, i.e.,
\begin{equation}
x_{ij}^ \al=x_{,ij}^ \al-\ch ijk x_k^ \al+ \cha  \bet \ga \al x_i^ \bet x_j^ \ga.
\end{equation}
The comma indicates ordinary partial derivatives.

In this implicit definition the \tit{second fundamental form} $(h_{ij})$ is taken
with respect to $-\n$.

The second equation is the \tit{Weingarten equation}
\begin{equation}
\n_i^ \al=h_i^k x_k^ \al,
\end{equation}
where we remember that $\n_i^ \al$ is a full tensor.

Finally, we have the \tit{Codazzi equation}
\begin{equation}
h_{ij;k}-h_{ik;j}=\riema \al \bet \ga \de\n^ \al x_i^ \bet x_j^ \ga x_k^ \de
\end{equation}
and the \tit{Gau{\ss} equation}
\begin{equation}
\riem ijkl= \{h_{ik}h_{jl}-h_{il}h_{jk}\} + \riema  \al \bet\ga \de x_i^ \al x_j^ \bet
x_k^ \ga x_l^ \de.
\end{equation}

When we consider hypersurfaces $M\su S^{n+1}$ to be embedded in $\R[n+2]$, we label the coordinates in $\R[n+2]$ as $(x^a)$, i.e., indices $a,b,c,...$ always run through $n+2$ values either from $1$ to $n+2$ or from $0$ to $n+1$.

Let us also state the definition of curvature functions of class $(K)$
\bd
A symmetric curvature function $F\in C^{2,\alpha}(\C_+)\ii C^0(\bar \C_+)$
positively homogeneous of degree $d_0>0$ is said to be of class $(K)$ if
\begin{equation}\lae{1.27}
F_i=\pd F{\kappa}i>0\q \text{in } \C_+,
\end{equation}
\begin{equation}\lae{1.28}
\fmo{\pa \C_+}=0,
\end{equation}
and
\begin{equation}\lae{1.29}
F^{ij,kl}\h_{ij}\h_{kl}\le F^{-1}(F^{ij}\h_{ij})^2-F^{ik}\tilde
h^{jl}\h_{ij}\h_{kl}\qq\A\,\h\in\msc S,
\end{equation}
or, equivalently, if we set $\hat F=\log F$,
\begin{equation}\lae{1.30}
\hat F^{ij,kl}\h_{ij}\h_{kl}\le -\hat F^{ik}\tilde
h^{jl}\h_{ij}\h_{kl}\qq\A\,\h\in\msc S,
\end{equation}
where $F$ is evaluated at $(h_{ij})$.
\ed

A detailed analysis of these curvature functions can be found in \cite[Section 1]{cg:indiana}. In this paper we actually do not need the full strength of inequality \re{1.30}. 

As we have shown in \ci[Lemma 1.3 and Remark 1.4]{cg97} a symmetric curvature
function $F\in C^2(\C_+)$ satisfies inequality \re{1.29} iff
\begin{equation}\lae{1.14}
F_i\kappa_i\le F_j\kappa_j\,,\q \text{for} \q \kappa_j\le \kappa_i,
\end{equation}
and
\begin{equation}\lae{1.20}
F_{ij}\x^i\x^j\le F^{-1}(F_i\x^i)^2-F_i\kappa_i^{-1}\abs{\x^i}^2\q \forall\,\x\in \R[n],
\end{equation}
where $F_i$, $F_{ij}$ are ordinary partial derivatives of $F$ in $\C_+$.

We only need the property \re{1.14}.

Let us finish this section with a simple yet useful observation.
\bl\lal{1.2}
Let $F\in (K)$ be homogeneous of degree $1$, then $F$ is concave.
\el

\bp
It suffices to prove that the right-hand side  of the inequality \re{1.20} is non-positive, if $F$ is homogeneous of degree $1$.

Using Schwarz's inequality we deduce
\begin{equation}
\begin{aligned}
F_i\xi^i&=\sum_i F_i^{\frac12}\ka_i^{\frac12}F_i^{\frac12}\ka_i^{-\frac12}\xi^i\\
&\le \big(\sum_iF_i\ka_i\big)^{\frac12} \big(\sum_iF_i\ka_i^{-1}\abs{\xi^i}^2\big)^{\frac12}=F^{\frac12} \big(\sum_iF_i\ka_i^{-1}\abs{\xi^i}^2\big)^{\frac12},
\end{aligned}
\end{equation}
hence the result.
\ep

\section{Polar sets}\las{2}

Let $M\su S^{n+1}$ be a connected, closed, immersed, strictly convex hypersurface given by an immersion
\begin{equation}
x:M_0\ra M\su S^{n+1},
\end{equation}
then $M$ is embedded, homeomorphic to $S^n$, contained in an open hemisphere and is the boundary of a convex body $\hat M\su S^{n+1}$, \cf \ci{docarmo}.

Considering $M$ as a codimension $2$ submanifold of $\R[n+2]$ such that
\begin{equation}\lae{2.2}
x_{ij}=-g_{ij}x-h_{ij}\tilde x,
\end{equation}
where $\tilde x\in T_x(\R[n+2])$ represents the exterior normal vector $\nu\in T_x(S^{n+1})$, we want to prove that the mapping
\begin{equation}\lae{2.3}
\tilde x:M_0\ra S^{n+1}
\end{equation}
is an embedding of a strictly convex, closed, connected hypersurface $\tilde M$. We call this mapping the \tit{Gau{\ss} map} of $M$.

First, we shall show that the Gau{\ss} map is injective. To prove this result we need the following lemma.
\bl\lal{2.1}
Let $M\su S^{n+1}$ be a closed, connected, strictly convex hypersurface and denote by $\hat M$ its (closed) convex body. Let $x\in M$ be fixed and $\tilde x$ be the corresponding outward normal vector, then
\begin{equation}\lae{2.4}
\spd y{\tilde x}\le 0\qq\A\, y\in \hat M
\end{equation}
and also strictly less than $0$ unless $y=x$.

The preceding inequality also characterizes the points in $\hat M$, namely, let $y\in S^{n+1}$ be such that
\begin{equation}\lae{2.5b}
\spd y{\tilde x}\le 0\qq\A\,x\in M,
\end{equation}
then $y\in \hat M$.
\el

\bp
\cq{\re{2.4}}\q   First, we note that $\hat M$ is contained in an open hemisphere $\mc H(x_0)$.

Let $y\in \tup{int}\, \hat M$ be arbitrary and let $z=z(t)$, $0\le t\le d$, be the unique minimizing geodesic in $S^{n+1}$ connecting $y$ and $x$ such that
\begin{equation}
z(0)=x\q\wed\q z(d)=y
\end{equation}
parametrized by arc length, and hence $0<d<\pi$.

Viewing $z$ as a curve in $\R[n+2]$ the geodesic equation has the form
\begin{equation}
\Ddot z\equiv \tfrac{D}{dt}\dot z =-z.
\end{equation}
If the coordinate system in $\R[n+2]$ is Euclidean, the covariant derivatives are just ordinary derivatives.

It is well-known that the geodesic $z$ is contained in $\hat M$ and that
\begin{equation}\lae{2.8}
\spd{\dot z(0)}{\tilde x}<0;
\end{equation}
notice that, after introducing geodesic polar coordinates in $S^{n+1}$ centered in $y$, we have
\begin{equation}
\spd{\dot z(0)}{\tilde x}=-\spd{\frac{\pa}{\pa r}}\nu
\end{equation}
and hence is strictly negative, \cf the remarks in \cite[Section 4]{cg96} after Theorem~4.6.

Thus, $\f(t)=\spd{z(t)}{\tilde x}$ satisfies the initial value problem
\begin{equation}\lae{2.9}
\Ddot \f=-\f,\q \f(0)=0,\q \dot\f(0)<0,
\end{equation}
and is therefore equal to
\begin{equation}
\f(t)=-\lam \sin t,\q \lam>0,
\end{equation}
i.e.,
\begin{equation}
\f(t)<0\qq\A\,0<t<\pi.
\end{equation}

Now, let $y\in M$, $y\ne x$, be arbitrary, and consider a sequence $z_k$ of geodesics parametrized in the interval $0\le t\le 1$, such that
\begin{equation}
z_k(0)=x\q\wed\q z_k(1)\ra y,
\end{equation}
where $z_k(1)\in \tup{int}\, \hat M$.

The geodesics $z_k$ converge to a geodesic $z$ connecting $x$ and $y$. If
\begin{equation}
\spd{\dot z(0)}{\tilde x}<0,
\end{equation}
then the previous arguments are valid yielding
\begin{equation}
\spd y{\tilde x}<0.
\end{equation}

On the other hand, the alternative
\begin{equation}
\spd y{\tilde x}=0
\end{equation}
leads to a contradiction, since then the geodesic $z$ would be part of the tangent space $T_x(M)$ which is impossible, \cf the considerations in  \cite{cg96} after the equation $(4.17)$.

\cvm
\cq{$y\in \hat M$}\q Suppose now that $y\in S^{n+1}$ satisfies \re{2.5b}, and assume by contradiction that $y\in \cpl \hat M$. Pick an arbitrary $\bar x_0\in \tup{int}\msp \hat M$, $\bar x_0\ne -y$, and let $z=z(t)$, $0\le t\le d$,  be the minimizing geodesic joining $\bar x_0$ and $y$ parameterized by arc length, such that $z(0)=\bar x_0$ and $z(d)=y$. The geodesic intersects $M$ in a unique point $x$, $x=z(t_1)$, $0<t_1<d$.

Define
\begin{equation}
\f(t)=\spd {z(t)}{\tilde x},
\end{equation}
then
\begin{equation}
\f(t_1)=0\q\wed\q \dot\f(t_1)>0,
\end{equation}
and hence
\begin{equation}
\f(t)=\lam \sin(t-t_1), \q \lam >0,
\end{equation}
and we  conclude
\begin{equation}
\f(t)>0\qq\A\, t_1<t<t_1+\pi
\end{equation}
contradicting the assumption  $\f(d)\le 0$.

Therefore we have proved $y\in\hat M$.
\ep

\bt\lat{2.2}
Let $x:M_0\ra M\su S^{n+1}$ be the embedding of a closed, connected, strictly convex hypersurface, then the Gau{\ss} map defined in \re {2.3} is injective, where we identify $\R[n+2]$ with its individual tangent spaces.
\et

\bp
We again assume $M$ to be a codimension $2$ submanifold in $\R[n+2]$. Suppose there would be two points $p_1\ne p_2$ in $M_0$ such that
\begin{equation}
\tilde x(p_1)=\tilde x(p_2),
\end{equation}
then the function
\begin{equation}
\f(y)=\spd y{\tilde x(p_1)}
\end{equation}
would vanish in the points $x(p_1)$ as well as $x(p_2)$ contrary to the results of \rl{2.1}.
\ep

\bl
As a submanifold of codimension $2$ $M$ satisfies the Wein\-gar\-ten equations
\begin{equation}\lae{2.18b}
\tilde x_i=h^k_ix_k
\end{equation}
for the normal $\tilde x$ and also
\begin{equation}
x_i=g^k_ix_k
\end{equation}
for the normal $x$.
\el

\bp
We only have to prove the non-trivial Weingarten equation.

First we infer from
\begin{equation}
\spd x{\tilde x}=0
\end{equation}
that
\begin{equation}
0=\spd{x_i}{\tilde x}+\spd x{\tilde x_i}=\spd x{\tilde x_i}.
\end{equation}
Furthermore, there holds
\begin{equation}\lae{2.22}
0=\spd{\tilde x}{\tilde x_i},
\end{equation}
since $\spd {\tilde x}{\tilde x}=1$. Hence, we deduce
\begin{equation}\lae{2.23b}
\tilde x_i=a^k_ix_k.
\end{equation}

Differentiating the relation $\spd{x_j}{\tilde x}=0$ covariantly we obtain
\begin{equation}\lae{2.24}
\spd {\tilde x_j}{x_i}=h_{ij}
\end{equation}
and we infer \re{2.18b} in view of \re{2.23b}.
\ep

We can now prove
\bt\lat{2.4}
Let $x:M_0\ra M\su S^{n+1}$ be a closed, connected, strictly convex hypersurface of class $C^m$, $m\ge 3$, then the Gau{\ss} map $\tilde x$ in \re{2.3} is the embedding of a closed, connected, strictly convex hypersurface $\tilde M\su S^{n+1}$ of class $C^{m-1}$. 

Viewing $\tilde M$ as a codimension $2$ submanifold in $\R[n+2]$, its Gaussian formula is
\begin{equation}\lae{2.25}
\tilde x_{ij}=-\tilde g_{ij}\tilde x-\tilde h_{ij} x,
\end{equation}
where $\tilde g_{ij}$, $\tilde h_{ij}$ are the metric and second fundamental form of the hypersurface $\tilde M\su  S^{n+1}$, and $x=x(\xi)$ is the embedding of $M$ which also represents the exterior normal vector of $\tilde M$. The second fundamental form $\tilde h_{ij}$ is defined with respect to the interior normal vector. 

The second fundamental forms of $M$, $\tilde M$ and the corresponding principal curvatures $\ka_i$, $\tilde \ka_i$ satisfy
\begin{equation}\lae{2.26b}
h_{ij}=\tilde h_{ij}=\spd{\tilde x_i}{x_j}
\end{equation}
and
\begin{equation}\lae{2.27b}
\tilde \ka_i=\ka_i^{-1}.
\end{equation}
\et

\bp
(i) From the Weingarten equation \re{2.18b} we infer
\begin{equation}\lae{2.28}
\tilde g_{ij}=\spd{\tilde x_i}{\tilde x_j}=h^k_ih_{kj}
\end{equation}
is positive definite, hence $\tilde x=\tilde x(\xi)$ is an embedding of a closed, connected hypersurface, where we also used \rt{2.2}.

\cvm
(ii) The pair $(x,\tilde x)$ satisfies
\begin{equation}\lae{2.30b}
\spd x{\tilde x}=0
\end{equation}
and we claim that $x$ is the exterior normal vector of $\tilde M$ in $\tilde x$, where as usual we identify the normal vector $\tilde \nu=(\tilde\nu^\al)\in T_{\tilde x}(S^{n+1})$ with its embedding in $T_{\tilde x}(\R[n+2])$.

Differentiating \re{2.30b} covariantly and using the fact that $\tilde x$ is a normal vector for $M$ we deduce
\begin{equation}\lae{2.31b}
0=\spd x{\tilde x_i},
\end{equation}
i.e., $\tilde x$ and $x$ span the normal space of the codimension $2$ submanifold $\tilde M$. 

Let us define the second fundamental form $\tilde h_{ij}$ of $\tilde M\su S^{n+1}$ with respect to the normal vector $\tilde \nu \in T_{\tilde x}(S^{n+1})$ corresponding to $x$, then the codimension $2$ Gaussian formula is exactly \re{2.25}.

Differentiating the Weingarten equation \re{2.18b} covariantly with respect to the metric $\tilde g_{ij}$ and indicating the covariant derivatives with respect to $\tilde g_{ij}$ by a semi-colon and those with respect to $g_{ij}$ simply by indices, we obtain
\begin{equation}
\tilde x_{;ij}=h^k_{i;j}x_k+h^k_ix_{;kj}
\end{equation}
and we deduce further
\begin{equation}\lae{2.33b}
\tilde h_{ij}=-\spd{\tilde x_{;ij}}x=-h^k_i\spd{x_{kj}}x=h^k_ig_{kj}=h_{ij}.
\end{equation}

On the other hand, we infer from \re{2.31b}
\begin{equation}
\tilde h_{ij}=-\spd{\tilde x_{;ij}}x=\spd{\tilde x_i}{x_j}
\end{equation}
which proves \re{2.26b}.

The last relation \re{2.27b} follows from \re{2.33b} and \re{2.28}. 

Finally, the normal vector $x$ must correspond to the exterior normal of $\tilde M$ in $T_{\tilde x}(S^{n+1})$, since $\tilde h_{ij}$ is positive definite.
\ep

We can also define a Gau{\ss} map from the strictly convex, connected, closed hypersurfaces $\tilde M\su S^{n+1}$ into $S^{n+1}$, and the preceding theorem shows that the two Gau{\ss} maps are inverse to each other,  i.e., if we start with a closed, strictly convex hypersurface $M\su S^{n+1}$, apply the Gau{\ss} map to obtain a  strictly convex hypersurface $\tilde M\su S^{n+1}$, and then apply the second Gau{\ss} map, then we return to $M$ with a pointwise equality.

Denoting the two Gau{\ss} maps simply by a tilde, this can be expressed in the form
\begin{equation}
x=\tilde{\tilde x},
\end{equation}
or, equivalently, in the form of a commutative diagram
\begin{equation}
\begin{diagram}
\node{M}
      \arrow[2]{e,t}{\tilde{}} \arrow{se,r}{\id}
   \node[2]{\tilde M} \arrow{sw,b}{\tilde{}}\\
\node[2]{M}
\end{diagram}
\end{equation}

Before we give an equivalent characterization of the images of the Gau{\ss} maps, let us show that the images of  strictly convex hypersurfaces by the Gau{\ss} maps are as smooth as the original hypersurfaces.

\bl\lal{2.5b}
Let $M\su S^{n+1}$ be a closed, connected, strictly convex hypersurface of class $C^{m,\al}$, $m\ge 3$, $0\le \al\le 1$ and let  $\tilde M\su S^{n+1}$ be its image under the Gau{\ss} map. Let $\tilde M\su \mc H(x_0)$ and express  $\tilde M$  as a graph in geodesic polar coordinates $(\rho,x^i)$ centered in $x_0$, $\tilde M=\fv{\graph \tilde u}{S^n}$, then $\tilde h_{ij}$, expressed in corresponding local coordinates $x^i$ of $S^n$, is of class $C^{m-2,\al}$.
\el

\bp
Notice that this is a non-trivial statement, since $\tilde M$ is only known to be of class $C^{m-1,\al}$. 

Let $(x^a)=(x^0,x^i)$ be Euclidean coordinates in $\R[n+2]$ and assume without loss of generality that $x_0=(1,0)$. Writing $x=(x^0,z)$, $z\in\R$, we have
\begin{equation}
\abs{x^0}^2=1-\abs z^2\qq\A\,x\in S^{n+1}, 
\end{equation}
i.e., after introducing Euclidean polar coordinates $(r,x^i)$ in $\R$, the hemisphere $\mc H(x_0)$ is given as the embedding
\begin{equation}\lae{2.42}
x=(x^0,r,x^i)=(\sqrt{1-r^2},r,x^i)
\end{equation}
and the lower hemisphere $\mc H(-x_0)$ by the embedding
\begin{equation}
x=(x^0,r,x^i)=(-\sqrt{1-r^2},r,x^i).
\end{equation}

The metric in $S^{n+1}\sminus \{x^0=0\}$ is then expressed as
\begin{equation}
d\bar s^2=\tfrac 1{1-r^2}dr^2+r^2\s_{ij}dx^idx^j,
\end{equation}
where $\s_{ij}$ is the metric of $S^n$. 

Defining $\rho$ by
\begin{equation}
d\rho=\frac1{\sqrt{1-r^2}}dr\q\wed\q \rho(0)=0
\end{equation}
will give us geodesic polar coordinates $(\rho,x^i)$ in $\mc H(x_0)$ centered in $x_0$. 

Now, assuming $\tilde M\su \mc H(x_0)$ implies $M\su \mc H(-x_0)$, in view of \rl{2.1}. Let $(\xi^i)$ be local coordinates for $M$ and express the Gau{\ss} map $\tilde x(\xi)$ in the coordinates in \re{2.42}
\begin{equation}
\tilde x(\xi)=(x^0(\xi),r(\xi),x^i(\xi)),
\end{equation}
then 
\begin{equation}
r(\xi)=u(x^i(\xi))\q\wed\q x^0(\xi)=\sqrt{1-u^2(x^i(\xi))},
\end{equation}
where $\tilde M$ has been written as a graph over $S^n$
\begin{equation}
\tilde M=\set{r=u(x^i)}{(x^i)\in S^n};
\end{equation}
notice that in geodesic polar coordinates we have $\tilde M=\graph \tilde u$ with 
\begin{equation}
\tilde u=\rho (u).
\end{equation}

In the coordinates $(\xi^i)$ the second fundamental form $\tilde h_{ij}$ is already known to be of class $C^{m-2,\al}$ because of the relation \re{2.26b}. Hence the lemma will be proved, if we can show that the transformation $(x^i(\xi))$ is a $C^{m-1,\al}$-diffeo\-morphism, i.e., we have to show that the Jacobian is invertible.

Now, the induced metric $\tilde g_{ij}$ can be expressed as
\begin{equation}
\begin{aligned}
\tilde g_{ij}&=\spd{\tilde x_i}{\tilde x_j}=x^0_ix^0_j+r_ir_j+r^2\s_{kl}x^k_ix^k_j\\
&=\tfrac1{1-r^2}r_ir_j+r^2\s_{kl}x^k_ix^l_j\\
&=\{\tfrac1{1-u^2}u_ku_l+u^2\s_{kl}\}x^k_ix^l_j,
\end{aligned}
\end{equation}
hence $(x^k_i)$ is invertible, since the left-hand side of this equation has this property.
\ep

\bt
Let $M\su S^{n+1}$ be a closed, connected, strictly convex hypersurface of class $C^{m,\al}$, $m\ge 2$, $0\le \al\le 1$, then $\tilde M\su N$, its image under the Gau{\ss} map is also of class $C^{m,\al}$. 
\et

\bp
(i) First, let us assume that $m\ge 3$ and $0\le \al\le 1$. The Gau{\ss} map is then of class $C^{m-1,\al}$, i.e., $\tilde M$ is of class $C^{m-1,\al}$. Here, we use the coordinates $(\xi^i)$ for $M$ also as coordinates for $\tilde M$. The metric $\tilde g_{ij}$ and the Christoffel symbols  of $\tilde M$ are then  of class $C^{m-2,\al}$ \resp $C^{m-3,\al}$, while the second fundamental form $\tilde h_{ij}$ is of class $C^{m-2,\al}$, in view of \re{2.26b}.

We may assume that $M\su \mc H(-x_0)$ and $\tilde M\su \mc H(x_0)$, where $x_0=(1,0)$. Using then geodesic polar coordinates $(\rho,\xi^i)$ centered in $x_0$, the metric in $S^{n+1}$ can be expressed in the form
\begin{equation}
d\bar s^2=d\rho^2+e^{2\psi(\rho)}\s_{ij}d\xi^id\xi^j,
\end{equation}
or, in conformal coordinates
\begin{equation}\lae{2.52}
d\bar s^2=e^{2\psi(\rho)}\{d\tau^2+\s_{ij}d\xi^id\xi^j\}.
\end{equation}

Writing $\tilde M$ as a graph in the coordinates $(\tau,\xi^i)$
\begin{equation}
\tilde M=\graph\fv u{S^n},
\end{equation}
the second fundamental form $h_{ij}$ of $\tilde M$ can be expressed as
\begin{equation}
e^{-\psi}v^{-1}h_{ij}=-u_{ij}-\cha 000\mspace{1mu}u_iu_j-\cha 0i0
\mspace{1mu}u_j-\cha 0j0\mspace{1mu}u_i-\cha ij0,
\end{equation}
where
\begin{equation}
v^2=1+\s^{ij}u_iu_j
\end{equation}
and where we note  that the second fundamental form $h_{ij}$  is of class $C^{m-2,\al}$, \cf \rl{2.5b}.

We want to replace the covariant derivatives $u_{ij}$ of $u$ with the covariant derivatives $u_{;ij}$ of $u$ with respect to the metric $\s_{ij}$  to deduce that $u_{;ij}$ is of class $C^{m-2,\al}$, and hence $u\in C^{m,\al}(S^n)$.

To achieve this we define a new metric $\hat g_{\al\bet}$ in the ambient space
\begin{equation}
\hat g_{\al\bet}=e^{-2\psi}\bar g_{\al\bet},
\end{equation}
where $\bar g_{\al\bet}$ is the metric in \re{2.52}. 
Let $\hat g_{ij}$, $\hat h_{ij}$ and $\hat \nu$ be the obvious geometric quantities of $\tilde M$ with respect to the new metric, then there holds
\begin{equation}
h_{ij} e^{-\psi}=\hat h_{ij}+\psi_\al\hat\nu^\al \hat g_{ij}
\end{equation}
as one easily checks.

On the other hand, $\hat h_{ij}$ can be expressed in terms of the Hessian $u_{;ij}$ of $u$ with respect to  the metric $\s_{ij}$, namely,
\begin{equation}
\hat h_{ij}=-u_{;ij}v^{-1},
\end{equation}
i.e.,
\begin{equation}
h_{ij}e^{-\psi}=-u_{;ij}v^{-1} +\psi_\al\hat\nu^\al(u_iu_j+\s_{ij}),
\end{equation}
hence, $u_{;ij}$ is of class $C^{m-2,\al}$.

\cvm
(ii) The case $m=2$ and $0\le \al\le 1$ follows by approximation and the uniform $C^{2,\al}$-estimates. Notice that the approximating second fundamental forms will converge in $C^0$.
\ep

\bd
(i) Let $M\su S^{n+1}$ be a closed, connected, strictly convex hypersurface, then we define its \tit{polar set} $M^*\su S^{n+1}$ by
\begin{equation}
M^*=\set{y\in S^{n+1}}{\sup_{x\in M}\spd xy=0},
\end{equation}
where the scalar product is the scalar  product in $\R[n+2]$ and $x$, $y$ are Euclidean coordinates. 

\cvm
(ii) Let $\hat M$ be the convex body of $M\su S^{n+1}$, then we define the polar of $\hat M$ by
\begin{equation}
\hat M^*=\set{y\in S^{n+1}}{\sup_{x\in \hat M}\spd xy\le 0}.
\end{equation}
\ed

\bt
The $M\su S^{n+1}$ be a closed, connected and strictly convex hypersurface, then
\begin{equation}\lae{2.62}
M^*=\tilde M
\end{equation}
and
\begin{equation}\lae{2.63}
\hat M^*=\hat{\tilde M}.
\end{equation}
\et

\bp
\cq{\re{2.62}}\q In view of \rl{2.1} there holds 
 \begin{equation}
\tilde M\su M^*.
\end{equation}

On the other hand, let $y\in M^*$ and $x\in M$ be such that
\begin{equation}\lae{2.53b}
\spd xy=0.
\end{equation}
Then we deduce, after introducing local coordinates in $M$,
\begin{equation}\lae{2.54b}
\spd{x_i}y=0
\end{equation}
and
\begin{equation}\lae{2.55b}
\spd{x_{ij}}y\le 0,
\end{equation}
where the derivatives are covariant derivatives with respect to the induced metric $g_{ij}$ of $M$ being viewed as a codimension $2$ submanifold.

Combining \re{2.53b} and \re{2.54b} we infer
\begin{equation}
y=\pm \tilde x,
\end{equation}
but because of \re{2.2} and \re{2.55b} we deduce $y=\tilde x$. 

\cvm
\cq{\re{2.63}}\q In view of \rl{2.1} we immediately deduce
\begin{equation}
\hat M^*\su \hat{\tilde M},
\end{equation}
hence we only have to prove the reverse inclusion.

Let $y\in \hat M$ and $\tilde x,\tilde z\in \tilde M$, $\tilde x\ne \tilde z$, be arbitrary and let $z=z(t)$ be the minimizing geodesic connecting $z(0)=\tilde x$ and $z(d)=\tilde z$ parametrized by arc length. 
Then it suffices to prove 
\begin{equation}
\f(t)=\spd y{z(t)}\le 0\qq\A\,0\le t\le d<\pi.
\end{equation}

Assume by contradiction that there exists $0<t_0<d$ such that
\begin{equation}
0<\f(t_0)=\sup\set{\f(t)}{0\le t\le d},
\end{equation}
then $\f$ solves the initial value problem
\begin{equation}
\Ddot\f=-\f,\q \f(t_0)>0,\q\dot\f(t_0)=0,
\end{equation}
and hence, it is equal to
\begin{equation}
\f(t)=\lam \cos (t-t_0)\q\lam>0,
\end{equation}
contradicting the relations $\f(0)\le 0$ and $\f(d)\le 0$, \cf \rl{2.1},  since there holds
\begin{equation}
0<t_0<\pih\q\vee \q 0<d-t_0<\pih.\qedhere
\end{equation}
\ep

An important corollary is
\bc
\tup{(i)} Let $M_i$, $i=1,2$, be connected, closed, strictly convex hypersurfaces in $S^{n+1}$, then
\begin{equation}
\hat M_1\su \hat M_2 \im \hat M_2^*\su \hat M_1^*.
\end{equation}

\tup{(ii)} Let $B_r(x_0)\su \mc H(x_0)$ be a geodesic ball of radius $0<r<\pih$, then its polar set is a closed geodesic ball centered in $-x_0$
\begin{equation}\lae{2.76}
B_r(x_0)^*=\bar B_{r^*}(-x_0),\qq 0<r^*=\f(r)<\pih,
\end{equation}
where $\f$ is continuous function.
\ec

\bp
We only need to prove \re{2.76}. But since the convex body of a geodesic sphere is the corresponding closed geodesic ball, it suffices to prove that the polar of a geodesic sphere $S_r(x_0)$ is a geodesic sphere $S_{r^*}(-x_0)$.

Let $\tilde M$ be the polar of $S_r(x_0)$, then we deduce from \re{2.26b}, that $\tilde M$ is totally umbilic and hence a geodesic sphere, \cf \rs{5} for details. This sphere must be centered in $-x_0$, since it is invariant under all $A\in O(n+2)$ having $x_0$ as a fixed point.
\ep

To conclude this section, we note that, with the help of the Gau{\ss} map,  the Minkowski type equation
\begin{equation}
\fv FM=f(\nu)
\end{equation}
in $S^{n+1}$ can be expressed in the form
\begin{equation}\lae{2.48}
\fv FM=f(\tilde x),
\end{equation}
where $f$ is supposed to be defined in $S^{n+1}$, or more precisely, in $T_x(\R[n+2])\equiv \R[n+2]$, the latter can be achieved by extending $f$ homogeneously of degree $0$.

Let $M^*$ be the polar set of $M$, $\tilde F$ the inverse of $F$, then the equation \re{2.48} is equivalent to
\begin{equation}\lae{2.49}
\fv {\tilde F}{M^*}=f^{-1}(\tilde x),
\end{equation}
where this time the right-hand side is looked at to be a function defined in the ambient space of $M^*$. Solving one equation is equivalent to solving the other.

\section{Curvature estimates}\las 3

We prove curvature estimates for the polar hypersurface $M^*$ satisfying the equation \re{2.49}. Since neither the result nor its proof relies on the fact that the underlying hypersurface is a polar hypersurface, we consider in this and in the following sections a strictly convex hypersurface $M$ satisfying the equation
\begin{equation}\lae{3.1}
\fv FM=f(x)\qq\A\, x\in M,
\end{equation}
where $0<f\in C^5(S^{n+1})$ and $F\in (K)$ of class $C^5$, and we shall prove that this problem has a solution, if  \ras{0.1} is satisfied. Since any positive power of $F$ is again of class $(K)$, we shall assume that $F$ is homogeneous of degree $1$ and hence concave, \cf \rl{1.2}. 
\bt\lat{3.1}
Let $M\in C^{4,\al}$ be a strictly convex hypersurface in $S^{n+1}$ satisfying the equation \re{3.1}, then its principal curvatures $\ka_i$ are uniformly bounded, i.e., there exist positive constants $c_1, c_2$ such that
\begin{equation}
0<c_1\le \ka_i\le c_2\qq\A\, 1\le i\le n,
\end{equation}
where the $c_i$ only depend on $F$ and $f$, which are supposed to satisfy the requirements mentioned above.
\et

\bp
It suffices to prove the upper estimate, since the lower estimate follows from the fact that $F$ is continuous in $\bar\C_+$ and vanishes on the boundary.

The second fundamental form $h_{ij}$ satisfies the equation
\begin{equation}\lae{3.3}
\begin{aligned}
-F^{kl}h^i_{j;kl}&=F^{kl}h_{kr}h^r_lh^i_j-Fh^{ki}h_{kj}+F^{kl,rs}h_{kl;j}h_{rs;m}g^{mi}\\[\cma]
&\hp{+}-f_{\al\bet}x^\al_kx^\bet_jg^{ki}+f_\al\nu^\al h^i_j+F\de^i_j-F^{kl}g_{kl}h^i_j,
\end{aligned}
\end{equation}
\cf the corresponding equation in \cite[equ. (5.4)]{cg:spaceform}, where an evolution problem is considered. The present situation can be recovered by setting $\dot\F=1$, $\Ddot \F=0$, $\tilde f=f$, $\F-\tilde f=0$ and $K_N=1$. 

We want to apply the maximum principle to obtain an a priori estimate for
\begin{equation}
\f=\sup\set{h_{ij}\h^i\h^j}{\norm\h=1}.
\end{equation}

Let $x_0\in M$ be a point where $\f$ attains its maximum. We then introduce Riemannian normal coordinates $\xi^i$ at $x_0$ such that at $x_0=x(\xi_0)$ we have
\begin{equation}
g_{ij}=\de_{ij},\q h_{ij}=\ka_i\de_{ij}\q\text{and}\q \f=h^n_n.
\end{equation}

Let $\h=(\h^i)$ be the contravariant vector field defined by
\begin{equation}
\h=(0,\ldots,1)
\end{equation}
in a neighbourhood of $\xi_0$ and set
\begin{equation}
\tilde\f=\frac{h_{ij}\h^i\h^j}{g_{ij}\h^i\h^j}.
\end{equation}
$\tilde\f$ is well defined in a neighbourhood of $\xi_0$.

Now $\tilde\f$ assumes its maximum at $\xi=\xi_0$. Moreover, at $\xi=\xi_0$ the covariant derivatives up to order two of $\tilde \f$ coincide with those of $h_n^n$, i.e., $\tilde\f$ satisfies the same differential equation at $\xi_0$ as $h^n_n$. For the sake of greater clarity let us therefore treat $h^n_n$ like a scalar and pretend that $\f$ is defined by
\begin{equation}
\f=h^n_n.
\end{equation}

Applying the maximum principle in $\xi_0$ we deduce
\begin{equation}
\begin{aligned}
0&\le F^{kl}h_{kr}h^r_lh^n_n-F\abs{h^n_n}^2+F^{kl,rs}h_{kl;n}h_{rs;m}g^{mn}\\[\cma]
&\hp{+}-f_{\al\bet}x^\al_kx^\bet_ng^{kn}+f_\al\nu^\al h^n_n+F-F^{kl}g_{kl}h^n_n,
\end{aligned}
\end{equation}
yielding
\begin{equation}
\begin{aligned}
0&\le F^{kl}h_{kr}h^r_lh^n_n-F\abs{h^n_n}^2+c_0(1+h^n_n)-F^{kl}g_{kl}h^n_n,
\end{aligned}
\end{equation}
where
\begin{equation}
c_0=c_0(\abs{f}_{{\vphantom{}}_{2,0}}).
\end{equation}

The function $F$ is of class $(K)$ and thus satisfies the estimate \re{1.14}. Let $\ka_1$ be the smallest principal curvature of $M$ in $x_0$, then
\begin{equation}\lae{3.12}
F=\sum_iF_i\ka_i\le nF_1\ka_1
\end{equation}
and hence
\begin{equation}
\begin{aligned}
F^{kl}h^r_kh_{rl}h^n_n-F\abs{h^n_n}^2&=F_1\ka_1(\ka_1-\ka_n)\ka_n+\sum_{i=2}^nF_i\ka_i(\ka_i-\ka_n)\ka_n\\[\cma]
&\le -\tfrac1n F (\ka_n-\ka_1)\ka_n.
\end{aligned}
\end{equation}

Now, if $\ka_n$ is supposed to be large in $x_0$, then
\begin{equation}
\ka_1\le \frac{\ka_n}2,
\end{equation}
because $F=f$ is bounded, hence $\ka_n=h^n_n$ is a priori bounded.
\ep

\section{Lower order bounds}\las 4

To derive the lower bounds we use the group invariance assumption. Let $M\su  S^{n+1}$ be a strictly convex, closed hypersurface and suppose that $M$ is invariant with respect to the group $G\su O(n+2)$
\begin{equation}
AM\su M\qq\A\, A\in G.
\end{equation}

Assume furthermore that a common fixed point $x_0$ of $G$ is an interior point of $\hat M$. The principal curvatures $\ka_i$ of $M$ are then also invariant with respect to $G$, i.e.,
\begin{equation}
\ka_i(x)=\ka_i(Ax)\qq\A\, x\in M,\; \A\, A\in G,
\end{equation}
as one easily checks.

Representing $M$ in  geodesic polar coordinates with center $x_0$ as a graph $u=u(\xi)$ over $S^n$, we conclude that the function $u$ is also invariant with respect to the induced isometry group in $S^n$, still denoted by $G$, i.e.,
\begin{equation}
u(\xi)=u(A\xi)\qq\A\, \xi\in S^n,\; \A\, A\in G.
\end{equation}

Since by assumption the induced group has no fixed points in $S^n$, $u$ is orthogonal to the first eigenfunctions of the Laplace operator in $S^n$, i.e.,
\begin{equation}
\int_{S^n}x^i u=0\qq\A\, 1\le i\le n+1,
\end{equation}
\cf \cite[Proposition 2.5]{guan:annals}.  Let us state this result as lemma.

\bl\lal{4.1}
Let $u\in C^0(S^n)$ be invariant with respect to the induced group $G$, then $u$ is orthogonal to the spherical harmonics of degree $1$.
\el

Now, we use stereographic projection $\pi$ to compare $M$ with a strictly convex hypersurface $\pi(M)\su \R$. Let $-x_0$ be the north pole of $S^{n+1}$ and assume that $\hat M$ is contained in the lower open hemisphere $\mc H(x_0)$
\begin{equation}
\hat M\su \mc H(x_0)
\end{equation}
such that $x_0\in \text{int\,} \hat M$, notice that by definition a convex body is always closed. The metric $\bar g_{\al\bet}$ of $S^{n+1}$ is then conformal to the Euclidean metric
\begin{equation}
\bar g_{\al\bet}=\frac1{(1+\tfrac14\abs{x}^2)^2}\de_{\al\bet},
\end{equation}
where $x=(x^\al)$ are Euclidean coordinates in $\R$.

The point $x_0\in S^{n+1}$ corresponds to the origin $0\in\R$ and, introducing Euclidean polar coordinates $(\rho,\xi^i)$, the metric in $S^{n+1}$ is expressed as
\begin{equation}
d\bar s^2=\frac1{(1+\tfrac14\rho^2)^2}\{d\rho^2+\rho^2\s_{ij}d\xi^id\xi^j\}.
\end{equation}

Comparing this expression with the representation of $\bar g_{\al\bet}$ in geodesic polar coordinates $(r,\xi^i)$ centered in $x_0$, namely,
\begin{equation}
d\bar s^2=dr^2+h(r)\s_{ij}d\xi^id\xi^j
\end{equation}
and observing that the radial geodesics in $S^{n+1}$ are mapped onto the radial geodesics in $\R$ we deduce that
\begin{equation}
r=\int_0^\rho\frac1{1+\tfrac14 t^2}=2 \arctan \tfrac\rho2.
\end{equation}

Finally, defining
\begin{equation}\lae{4.10}
\tau=\log \rho,
\end{equation}
we can express the metric in $S^{n+1}$ as
\begin{equation}\lae{4.11}
d\bar s^2=\frac{\rho^2}{(1+\tfrac14\rho^2)^2}\{d\tau^2+\s_{ij}d\xi^id\xi^j\}.
\end{equation}

Writing $M$ in these coordinates as a graph over $S^n$
\begin{equation}
M=\graph \fv u{S^n}
\end{equation}
$u$ is still invariant with respect to the induced group, and $\graph u$ also represents $\pi(M)$.

Let
\begin{equation}
\psi=-\log (1+\tfrac14\rho^2),
\end{equation}
such that
\begin{equation}
\bar g_{\al\bet}=e^{2\psi}\hat g_{\al\bet},
\end{equation}
where $(\hat g_{\al\bet})$ is the Euclidean metric, then the respective second fundamental forms $h_{ij}$ and $\hat h_{ij}$ are related by
\begin{equation}\lae{4.15}
e^\psi h^j_i=\hat h^j_i+\psi_\al\nu^\al\de^j_i,
\end{equation}
where $\nu$ is the exterior normal of $\pi(M)$ and
\begin{equation}
\psi_\al\nu^\al=\psi_0\nu^0=-\tfrac12\frac\rho{1+\tfrac14\rho^2}v^{-1},
\end{equation}
with
\begin{equation}
v^2=1+\s^{ij}u_iu_j.
\end{equation}

Thus, $\hat h_{ij}$ is also positive definite and invariant with respect to the induced group, as is the metric
\begin{equation}
\hat g_{ij}=e^{2u}\{u_iu_j+\s_{ij}\}.
\end{equation}

Moreover, since $\hat M$ is contained in the lower hemisphere, we have
\begin{equation}
0\le r\le \pih
\end{equation}
and hence
\begin{equation}
\rho\le 2\tan \frac\pi4=2.
\end{equation}

Thus, if the principal curvatures of $M$ are bounded by
\begin{equation}\lae{4.21}
0<k_1\le \ka_i\le k_2,
\end{equation}
then those of $\pi(M)$ are bounded by
\begin{equation}
0<\hat k_1\le \hat \ka_i\le \hat k_2,
\end{equation}
where
\begin{equation}
\hat k_j=\hat k_j(k_1,k_2),\qq j=1,2.
\end{equation}

Now we can prove that the convex body of $\pi(M)$ contains a Euclidean ball $B_{\rho_0}(0)$  and therefore $\hat M$ a geodesic ball $B_{r_0}(x_0)$.

\bl\lal{4.2}
Assume $x_0\in \tup{int\,}\hat M$, $\hat M\su \mc H(x_0)$, that $M$ is invariant with respect to the group $G$ and the principal curvatures satisfy the estimate \re{4.21}. Then there exists $0<r_0=r_0(k_1,k_2)$ such that the geodesic ball
\begin{equation}
B_{r_0}(x_0)\Su \tup{int\,}\hat M.
\end{equation}
\el

\bp
We shall prove that there exists a Euclidean ball of radius $0<\rho_0=\rho_0(\hat k_1,\hat k_2)$ such that
\begin{equation}\lae{4.25}
B_{\rho_0}(0)\Su \tup{int\,}\widehat{\pi(M)}.
\end{equation}

Let $\hat K$ be the Gaussian curvature of $\pi(M)=\graph u$, then $u, \,\hat g_{ij}$ and
\begin{equation}
\hat K=\hat K(u,\xi)
\end{equation}
are invariant functions in $S^n$ with respect to the induced group $G$, and hence orthogonal to the spherical harmonics of degree $1$, \cf \rl{4.1}. Hence the Steiner point $p$ of $\pi(M)$ coincides with the origin, since in Euclidean coordinates
\begin{equation}
\begin{aligned}
p^i&=\frac1{\abs{S^n}}\int_{\pi(M)}x^i\hat K\\[\cma]
&=\frac1{\abs{S^n}} \int_{S^n}x^i e^u\hat Kv=0.
\end{aligned}
\end{equation}

The relation \re{4.25} is then proved in \cite{cheng-yau:minkowski}. A similar, more general, result was later proved in \cite{schneider:balls}.
\ep

Let $M^*$ be the polar hypersurface of $M$, which is then also invariant with respect to the group $G$, since
\begin{equation}
0=\spd x{\tilde x}=\spd{Ax}{A\tilde x}\qq\A\, x\in M,\; \A\, A\in G.
\end{equation}
Then we shall prove

\bl
Let $M\su \mc H(x_0)$ be a strictly convex hypersurface satisfying the assumptions of \rl{4.2}. Then the polar convex body $\hat M^*$ of $\hat M$ is contained in $\mc H(-x_0)$ and there exist radii $0<r_1^*<r_0^*<\pih$ such that
\begin{equation}
B_{r_1^*}(-x_0)\Su \tup{int\,}\hat M^*\Su B_{r_0^*}(-x_0)\Su \mc H(-x_0).
\end{equation}
\el

\bp
Since $\hat M$ is compact there exists a geodesic ball $B_{r_1}(x_0)$ such that
\begin{equation}
\hat M \su B_{r_1}(x_0)\Su\mc H(x_0).
\end{equation}

Moreover, due to \rl{4.2}, there exists a geodesic ball $B_{r_0}(x_0)$ such that
\begin{equation}
B_{r_0}(x_0)\Su \tup{int\,}\hat M,
\end{equation}
hence we conclude
\begin{equation}
B_{r_1^*}(-x_0)=\tup{int}\msp  B^*_{r_1}(x_0)\Su \tup{int\,}\hat M^*\Su B^*_{r_0}(x_0)= \bar B_{r_0^*}(-x_0)\Su \mc H(-x_0). 
\end{equation}
\ep

Combining the two lemmata, and having in mind that both $M$ and $M^*$ are invariant with respect to $G$, so that \rl{4.2} can be applied to $M$ as well as $M^*$, we obtain 

\bt\lat{4.4}
Let $M\su \mc H(x_0)$ be a strictly convex hypersurface, invariant with respect to the group $G$ and assume that $x_0\in \tup{int\,}\hat M$ and that the principal curvatures $\ka _i$ satisfy the estimate \re{4.21}, then there exist radii $0<r_0<r_1<\pih$, depending only on the constants $k_j$, $j=1,2$, in \re{4.21} such that
\begin{equation}
B_{r_0}(x_0)\Su\tup{int\,}\hat M\Su B_{r_1}(x_0).
\end{equation}

The dual relation then also holds for $\hat M^*$, namely,
\begin{equation}
B_{r_1^*}(-x_0)\Su \tup{int\,}\hat M^*\Su  B_{r_0^*}(-x_0),
\end{equation}
where
\begin{equation}
\bar B_{r_i^*}(-x_0)=B^*_{r_i}(x_0),\qq i=0,1,
\end{equation}
and $0<r_1^*<r_0^*<\pih$.
\et

\section{A uniqueness result}\las 5

In this section we shall show that a strictly convex solution $M\su S^{n+1}$ of the equation
\begin{equation}\lae{5.1}
F=c\equiv\const>0,
\end{equation}
where $F$ is an arbitrary curvature function, homogeneous of degree $1$ and concave, is a geodesic sphere; notice that a curvature function is always supposed to be symmetric and monotone.

\bt\lat{5.1}
Let $M\su S^{n+1}$ be a closed strictly convex solution of \re{5.1}, then $M$ is a geodesic sphere. Assuming that $M$ is invariant with respect to the group $G$ and contained in $\mc H(x_0)$, where $x_0$ is a fixed point of $G$, then $M$ has to be a geodesic sphere with center in $x_0$.
\et

\bp
Assume without loss of generality that
\begin{equation}
F(1,\ldots,1)=n
\end{equation}
and consider the equation \re{3.3} for the second fundamental form. At a point $\bar x\in M$, where
\begin{equation}
\sup_M \max_i \ka_i=\ka_n=h^n_n
\end{equation}
is attained, the maximum principle yields, compare the proof of \rt{3.1},
\begin{equation}
\begin{aligned}
0&\le F^{kl}h_{kr}h^r_lh^n_n-F\abs{h^n_n}^2+F-F^{kl}g_{kl}h^n_n\\[\cma]
&=\sum_i F_i\ka_i(\ka_i-\ka_n)\ka_n+\sum_i F_i(\ka_i-\ka_n)\le 0.
\end{aligned}
\end{equation}
Hence $\bar x$ must be an umbilic and
\begin{equation}
c=F(\ka,\ldots,\ka)=\ka n,
\end{equation}
i.e.,
\begin{equation}
\sup_M\max_i\ka_i=\tfrac cn.
\end{equation}

But then all other points have to be umbilics too, since
\begin{equation}
c=F(\ka_i)\le F(\tfrac cn,\dots,\tfrac cn)=c.
\end{equation}

Now, any convex umbilic hypersurface $M$ of $S^{n+1}$ has to be a geodesic sphere, as can be most easily seen by choosing a point $y_0\in \hat M$ and using stereographic projection as in \rs 4. From equation \re{4.15} we then deduce that the projected hypersurface in Euclidean space is also umbilic and hence a sphere, \cf \cite[Vol. IV, p. 11]{spivak}.

If $M$ is invariant with respect to $G$ and contained in $\mc H(x_0)$, then its polar $M^*$ is also a strictly convex umbilic hypersurface such that its convex body contains a geodesic ball centered in $-x_0$ in its interior
\begin{equation}
B_{r^*_0}(-x_0)\Su \tup{int\,}\hat M^*,
\end{equation}
since  $\hat M$ is contained in a geodesic ball $B_{r_1}(x_0)\Su\mc H(x_0)$, in view of the compactness of $\hat M$ and the assumption $\hat M\su\mc H(x_0)$. Now for our purpose $M^*$ is as good as $M$, thus let us assume without loss of generality that $B_{r_0}(x_0)\Su \hat M$ and let us discard the assumption $\hat M\Su \mc H(x_0)$, since the corresponding result isn't known yet for $\hat M^*$.

Looking again at the stereographic projection $\pi(M)$ of $M$, where $x_0$ is now the south pole, i.e., $\pi(x_0)=0$, we still deduce that $\pi(M)$ is umbilic and hence a sphere, which now is invariant with respect to the group $G$. But as in the proof of \rl{4.2} we can then show that the Steiner point of $\pi(M)$ is the origin, and hence the origin must be the center of the sphere as one easily checks.

We then conclude that $M$ is a \cq{geodesic} sphere centered in $x_0$ by the properties of the stereographic projection. Using now the convexity of $M$ and the fact that $x_0$ is supposed to be part of $\hat M$, we obtain the final result that $\hat M\su \mc H(x_0)$ and that $M$ is a geodesic sphere centered in $x_0$.
\ep

\section{Existence of a solution}\las 6

The existence is proved by a continuity method using Smale's infinite dimensional version of Sard's theorem  \cite{smale:sard}. Writing the strictly convex hypersurfaces as graphs over $S^n$ it is convenient to express the differential operator
\begin{equation}
F=F(h_{ij})=\fv FM
\end{equation}
in terms of the standard Levi-Civit\`a connection in $S^n$.

Let $x_0\in S^n$ be a fixed point for the group $G$ and $\mc H(x_0)$ the corresponding hemisphere. Introducing geodesic polar coordinates centered in $x_0$, the metric in $\mc H(x_0)\sminus \{x_0\}$ can be expressed as
\begin{equation}
d\bar s^2=dr^2+e^{2\psi}\s_{ij}d\xi^id\xi^j,
\end{equation}
where $\psi=\psi(r)$, or in conformal coordinates
\begin{equation}
d\bar s^2=e^{2\psi}\{d\tau^2+\s_{ij}d\xi^id\xi^j\},
\end{equation}
where
\begin{equation}\lae{6.4}
\tau=\int_{\bar r}^re^{-\psi(t)},\qq 0<\bar r\le r<\pih,
\end{equation}
and $\bar r$ very small. Since all hypersurfaces we are concerned with lie in a region
\begin{equation}
\mc H(x_0,r_0,r_1)=\set{x\in \mc H(x_0)}{0<r_0\le r\le r_1<\pih},
\end{equation}
in view of \rt{4.4}, choosing $\bar r<r_0$ ensures that we do not have to worry about a possible coordinate singularity and still have a positive $\tau$-coordinate.

Let $M\su \mc H(x_0,r_0,r_1)$ be a strictly convex hypersurface, then, writing $M$ as a graph
\begin{equation}
M=\graph u=\set{(\tau,\xi)}{\tau=u(\xi), \xi\in S^n},
\end{equation}
the induced metric and the second fundamental form of $M$ are given by
\begin{equation}
g_{ij}=e^{2\psi}\{u_iu_j+\s_{ij}\}
\end{equation}
and
\begin{equation}
h_{ij}e^{-\psi}=\tilde h_{ij}+\psi_\bet\tilde\nu^\bet \tilde g_{ij},
\end{equation}
where the symbols with the tilde refer to the geometric quantities of $M$, when $M$ is considered to be embedded in the ambient space with metric
\begin{equation}
d\tilde s^2=d\tau^2+\s_{ij}d\xi^id\xi^j.
\end{equation}

$\tilde h_{ij}$ is then given by the relation
\begin{equation}
v^{-1}\tilde h_{ij}=-u_{ij}=-v^{-2}u_{;ij},
\end{equation}
where $u_{ij}$ is the Hessian of $u$ with respect to the induced metric $\tilde g_{ij}$ and $u_{;ij}$ is the Hessian of $u$ with respect to the standard metric $\s_{ij}$ of $S^n$. The term $v$ is defined by
\begin{equation}
v^2=1+\s^{ij}u_iu_j.
\end{equation}

Writing $u_{ij}$ instead of $u_{;ij}$ in the following, we see that
\begin{equation}\lae{6.12}
h_{ij}e^{-\psi}=-v^{-1}u_{ij}+v^{-1}\dot\psi\tilde g_{ij},
\end{equation}
where
\begin{equation}
\dot\psi=\df\psi\tau.
\end{equation}

If $M$ is invariant under $G$, then the function $u$ is also  invariant under the group action. Let $A^k(\xi)$ be the local representation of $A\xi$ and $(A^k_i)$ its derivative, then the covariant derivatives of $u$ satisfy
\begin{equation}\lae{6.14}
u_i(\xi)=u_k(A\xi)A^k_i,
\end{equation}
and
\begin{equation}\lae{6.15}
u_{ij}(\xi)=u_{kl}(A\xi)A^k_iA^l_j,
\end{equation}
notice that $A^k_{i;j}=0$.

\bd
A tensor field $\f$ in $T(S^n)$ is called invariant with respect to $G$, if it satisfies transformation relations according to \re{6.14}, \re{6.15}, where the contravariant indices transform like
\begin{equation}
\f^i(\xi)=\f^k(A\xi)A^i_k\qq\A\,A\in G,
\end{equation}
and there holds
\begin{equation}
A^i_kA^k_j=\de^i_j. 
\end{equation}
These transformation rules hold for invariant tensor fields of arbitrary order.
\ed

The metric $\s_{ij}$ of $S^n$ is of course invariant by the very definition of an isometry. Hence we conclude from \re{6.12} that the second fundamental form is also invariant and consequently also the tensor
\begin{equation}
F^{ij}=\pde F{h_{ij}}.
\end{equation}

Now consider the Banach spaces $E_1$, $E_2$ defined by
\begin{equation}
E_1=\set{u\in H^{5,p}(S^n)}{u\,\tup{invariant}}
\end{equation}
and
\begin{equation}
E_2=\set{w\in H^{3,p}(S^n)}{w\,\tup{invariant}}
\end{equation}
for some fixed $n<p<\un$, such that $H^{m,p}(S^n)\hra C^{m-1,\al}(S^n)$. 

Let $\Om\su E_1$ be an open bounded set such that $M(u)=\graph u$ is uniformly strictly convex, contained in $\mc H(x_0,r_0,r_1)$, such that $x_0$ is in interior point of $\widehat{M(u)}$ for all $u\in\Om$. We then define
\begin{equation}
\F:\Om\ra E_2
\end{equation}
by
\begin{equation}
\F(u)=F(h_{ij})-f(u,\xi)
\end{equation}
expressing a position vector $x\in \mc H(x_0)$ by $x=(\tau,\xi)$. 

All possible solutions of $\F=0$ are strictly contained in $\Om$, if $\Om$ is specified by the requirements
\begin{equation}
0<\tau_0=\tau(r_0)<u<\tau_1=\tau(r_1),
\end{equation}
\begin{equation}
x_0\in\tup{int\,}\widehat{M(u)},
\end{equation}
and
\begin{equation}
0<\e_0<\ka_i<\bar\ka,
\end{equation}
where $\ka_i$ are the principal curvatures of $\graph u$, in view of the a priori estimates in \rs 3 and \rs 4.

\bl
$\F$ is a proper nonlinear Fredholm operator of index zero.
\el

\bp
$F$ and hence $\F$ are uniformly elliptic in $\Om$. The properness is due to the a priori estimates in \rs 3 and \rs 4, the Evans-Krylov and Calder\`on-Zygmund  estimates and our assumption that $F$ and $f$ are of class $C^5$.

If the Banach spaces $E_i$ would have been defined without the symmetry requirement, the other properties of $\F$ would have been well known. Let $L$ be the derivative of $\F$, then $L$ is an elliptic linear partial differential operator of second order 
\begin{equation}
Lu=-F^{ij}u_{ij}+b^iu_i+cu
\end{equation}
and the lemma will be proved, if we can show that the operator
\begin{equation}
-F^{ij}u_{ij}+\lam u,\qq\lam>0,
\end{equation}
is surjective, i.e., for arbitrary $w\in E_2$ there exists $u\in E_1$ such that
\begin{equation}\lae{6.27}
-F^{ij}u_{ij}+\lam u=w.
\end{equation}

It is well known that there exists a function $u\in H^{5,p}(S^n)$ that solves the preceding equation, and we shall show $u$ is invariant, if $w$ is.

Let $A\in G$, then we claim that $\tilde u=u\circ A$ also satisfies \re{6.27}, which would yield
\begin{equation}
\tilde u=u
\end{equation}
because of the uniqueness.

Now, differentiating $\tilde u=u\circ A$ we obtain
\begin{equation}
\tilde u_{ij}(\xi)=u_{kl}(A\xi)A^k_iA^l_j
\end{equation}
and we infer
\begin{equation}
\begin{aligned}
-F^{ij}\tilde u_{ij}=-F^{ij}A^k_iA^l_j u_{kl} = -F^{kl}u_{kl},
\end{aligned}
\end{equation}
since $F^{ij}$ is invariant.
\ep

Recall that $w\in E_2$ is said to be a \tit{regular} value for $\F$, if either $w\notin  R(\F)$, or if for any $u\in \F^{-1}(w)$ $D\F(u)$ is surjective.

Smale \cite{smale:sard} proved that for separable Banach spaces $E_i$ and proper Fredholm maps $\F$ the set of regular values in $E_2$ is open and dense, if $\F$ is of class $C^k$ such that
\begin{equation}
k>\max(\tup{ind\,}\F,0).
\end{equation}
All requirements are satisfied in the present situation.

Next we want to use the uniqueness result in \rt{5.1}. Let $c>0$ be a constant such that
\begin{equation}\lae{6.32}
c<\inf_{S^n} f
\end{equation}
and let $u_0\equiv \const$ be such that the geodesic sphere $M_0=\graph u_0$ satisfies
\begin{equation}
\fv F{M_0}=c.
\end{equation}
We assume furthermore that the constants $r_0,r_1$ and $\e_0$, $\bar\ka$ are chosen such that all possible solutions of
\begin{equation}
F=tf+(1-t)c,\qq 0\le t\le 1,
\end{equation}
in $\mc H(x_0)$ satisfy the corresponding estimates.

The requirement \re{6.32} is not essential, it will only simplify some of the following arguments.

Let $0<\de$ be small and define
\begin{equation}
\Lam: \Om\ti [-\de,1+\de]\ra E_2
\end{equation}
 by
 \begin{equation}
\Lam (u,t)=F(h_{ij})-(tf+(1-t)c).
\end{equation}
Then $\Lam$ is also a proper Fredholm operator such that $\tup{ind\,}\Lam(\cdot,t)=0$ for fixed $t$, and, if $w\in E_2$ is a regular value for $\Lam$, then
\begin{equation}\lae{6.37}
\tup{ind\,}\Lam=1\qq\A\, (u,t)\in\Lam^{-1}(w).
\end{equation}

Recall that
\begin{equation}
\tup{ind\,}\Lam =\dim N(D\Lam)-\dim \tup{Coker\,} (D\Lam).
\end{equation}

The relation \re{6.37} will be proved in \rl{6.5} below.
\bt
Let $0<f\in C^5(S^n)$ be invariant under $G$, then for any $F\in (K)$ of class $C^5$, there exists a strictly convex invariant hypersurface $M\su \mc H(x_0)$ satisfying
\begin{equation}
\fv FM=f.
\end{equation}
\et

\bp
Consider the Fredholm map $\Lam =\Lam(u,t)$. The theorem will be proved, if we can show that there exists $u\in\Om$ such that
\begin{equation}
\Lam(u,1)=0.
\end{equation}

On the other hand, there exists a unique solution of the equation
\begin{equation}\lae{6.47}
\Lam(u,0)=0,
\end{equation}
namely, $u=u_0$, the geodesic sphere. In the lemma below we shall show that $u_0$ is also a regular point for $\Lam(\cdot,0)$, or equivalently, $(u_0,0)$ a regular point for $\Lam$.

Without loss of generality we may assume $0\notin R(\Lam(\cdot,1))$, for otherwise we have nothing to prove, and thus, $0$ is also regular value for $\Lam(\cdot,1)$.

Let $\e>0$ be small, then there exists a
\begin{equation}
w_\e\in B_\e(0)\su E_2,
\end{equation}
such that
\begin{equation}
tf+(1-t)c+w_\e>0\qq\A\,-\de\le t\le 1+\de,
\end{equation}
$w_\e\in R(\Lam(\cdot,0))$, and such that  $w_\e$ is a regular value for $\Lam(\cdot,0)$, $\Lam(\cdot,1)$ and $\Lam$.

Consider
\begin{equation}
\C_\e=\Lam^{-1}(w_\e)\ii (E_1\ti (-\de,1+\de)),
\end{equation}
then $\C_\e\ne\eS$ and $\C_\e$ is a $1$-dimensional submanifold without boundary.

The intersection
\begin{equation}
\tilde\C_\e=\C_\e\ii (E_1\ti [0,1])
\end{equation}
is then compact, since $\Lam$ is proper, and it consists of finitely many closed curves or segments.

We want to prove that there is $u_\e\in\Om$ such that $(u_\e,1)\in \tilde \C_\e$. Suppose this were not the case, then consider a point $(\bar u_\e,0)\in\tilde\C_\e$. Such points exist by assumption. Moreover, the $1$-dimensional connected submanifold $M_\e\su \C_\e$ containing $(\bar u_\e,0)$ can be expressed near $(\bar u_\e,0)$ by
\begin{equation}\lae{6.52}
M_\e=\set{(\f(t),t)}{-\de<t<\de},
\end{equation}
where $\f\in C^1$, $\f(0)=\bar u_\e$, and
\begin{equation}
\Lam(\f(t),t)=w_\e,
\end{equation}
since by assumption $D_1\Lam(\bar u_\e,0)$ is an isomorphism and the implicit function theorem can be applied.

Let $\tilde M_\e=M_\e\ii\tilde \C_\e$, then $\tilde M_\e$ isn't closed because of \re{6.52}, and hence has two endpoints, see \cite[Appendix]{milnor}. One of them is $(\bar u_\e,0)$ and the other also belongs to $\Lam(\cdot,0)^{-1}(w_\e)$ and can therefore be expressed as
\begin{equation}
(\tilde u_\e,0),
\end{equation}
where $\tilde u_\e\ne \bar u_\e$ because of the implicit function theorem.

Hence we have proved that the assumption
\begin{equation}
\Lam(\cdot,1)^{-1}(w_\e)=\eS
\end{equation}
implies
\begin{equation}
\#\Lam(\cdot,0)^{-1}(w_\e) \;\tup{is even}.
\end{equation}

However, we shall show that $\Lam(\cdot,0)^{-1}(w_\e)$ contains only one point, if $\e$ is small.

Indeed, let $\bar u_\e\in\Lam(\cdot,0)^{-1}(w_\e)$, then the $\bar u_\e$ converge to the unique solution $u_0$ of \re{6.47}. Thus, if $\e$ is small, all $\bar u_\e$ are contained in an open ball
\begin{equation}
B_\rho(u_0)\su\Om,
\end{equation}
where $\F=\Lam(\cdot,0)$ is a diffeomorphism due to \rl{6.4}, hence there exists just one solution of the equation 
\begin{equation}
\Lam(\bar u_\e,0)=w_\e.
\end{equation}

Thus we have proved that there exists a sequence
\begin{equation}
u_\e\in \Lam(\cdot,1)^{-1}(w_\e),
\end{equation}
if $\e$ tends to zero. A subsequence will then converge to a solution $u$ of
\begin{equation}
\Lam(u,1)=0.\qedhere
\end{equation}
\ep

It remains to prove two lemmata:
\bl\lal{6.4}
$u_0$ is a regular point for $\Lam(\cdot,0)$.
\el

\bp
Let $\f\in E_1$ be arbitrary and $\e>0$ so small that
\begin{equation}
u=u_0+\e\f\in\Om.
\end{equation}
Then we have to calculate
\begin{equation}
\df{}\e\fv{\Lam(u,0)}{\e=0}=\df{}\e\{F(h_{ij})-c\}.
\end{equation}

Now,
\begin{equation}
\df{}\e F(h_{ij})=F^i_j\dot h^j_i 
\end{equation}
and
\begin{equation}
h^j_i=-v^{-1}e^{-\psi}\tilde g^{jk}u_{ik}+\dot\psi e^{-\psi}\de^j_i,
\end{equation}
in view of \re{6.12}, where
\begin{equation}
\tilde g^{jk}=\s^{jk}-\frac{u^ju^k}{v^2}.
\end{equation}

Evaluating the resulting expressions at $\e=0$ we conclude
\begin{equation}
\dot h^j_i=-e^{-\psi}\f^j_i+\{\Ddot\psi-\abs{\dot\psi}^2\}e^{-\psi}\de^j_i\f,
\end{equation}
hence,
\begin{equation}
\df{}\e F(h_{ij})=e^{-\psi}\{-\D\f-n(\abs{\dot\psi}^2-\Ddot\psi)\f\},
\end{equation}
where the Laplace operator is taken with respect to the metric in $S^n$ and $e^{-\psi}$ is a constant.

Looking at the equations \re{4.10}, \re{4.11} we deduce that $\psi$ can be expressed as a function of $\tau$ as
\begin{equation}
\psi=\log\rho-\log(1+\tfrac14\rho^2),\q\rho=e^{\tau+\tau_0},
\end{equation}
where $\tau_0$ is an integration constant depending on the value of $\bar r$ in \re{6.4}, yielding
\begin{equation}
\abs{\dot\psi}^2-\Ddot\psi=1.
\end{equation}

Thus $\f\in N(D_1\Lam(u_0,0))$ satisfies
\begin{equation}
-\D\f-n\f=0
\end{equation}
and is therefore a spherical harmonic of degree $1$ or identically zero. But the $G$-invariant functions are orthogonal to the spherical harmonics of degree $1$, hence $D_1\Lam(u_0,0)$ is an isomorphism.
\ep
\bl\lal{6.5}
Let $\Lam$ be defined as above, then
\begin{equation}
\tup{ind\,}\Lam=1.
\end{equation}
\el

\bp
Let $(u_0,t_0)\in E_1\times [-\de,1+\de]$ be fixed, where we may assume that $t_0=1$, since $\Lam$ is continuous. We distinguish two cases:

\cvm
\tit{Case} 1: \q\; $(f-c)\in R(D\F(u_0))$

\cvm
Notice that $\Lam$ can be extended as a class $C^2$-function  to $E_1\ti \R[]$. We have
\begin{equation}
D\Lam=(D_1\Lam,-(f-c)),
\end{equation}
where all derivatives are evaluated at $(u_0,1)$ \resp $u_0$.  Then we deduce
\begin{equation}
\dim N(D\Lam)=\dim N(D_1\Lam)+1=\dim N(D\F)+1,
\end{equation}
for let 
\begin{equation}
D_1\Lam u_1=f-c,
\end{equation}
then
\begin{equation}
N(D\Lam)=N(D\F)\times \{0\}\oplus \langle{(u_1,1)}\rangle
\end{equation}
as one easily checks, and of course there holds
\begin{equation}
R(D\Lam)=R(D\F).
\end{equation}

\cvm
\tit{Case} 2:\q\; $(f-c)\notin R(D\F(u_0))$

\cvm
In this case
\begin{equation}
R(D\Lam)=R(D_1\Lam)\oplus \langle{(f-c)}\rangle
\end{equation}
and
\begin{equation}
N(D\Lam)=N(D_1\Lam)\times \{0\},
\end{equation}
hence
\begin{equation}
\tup{ind\,}\Lam=\tup{ind\,} \F+1=1
\end{equation}
in both cases.
\ep

\section{Proof of \rt{0.4}}\las 7 

The barrier condition for the original pair $(F,f)$ in $\mc H(-x_0)$  immediately translates to a barrier condition for $(\tilde F,f^{-1})$ in $\mc H(x_0)$. Following the stipulations in \rr{0.5}, we again assume that we consider the problem
\begin{equation}\lae{7.1}
\fv FM=f(x)\qq\A\,x\in M,
\end{equation}
where $F\in (K)$ and $M_1$ \resp $M_2$ are lower \resp upper barriers for $(F,f)$ bounding a connected open set $\Om\su \mc H(x_0)$.

We want to apply an old result, \cite[Theorem 0.4]{cg97}, in which we showed that the problem \re{7.1} has a strictly convex solution $M\su \bar \Om$ of class $C^{6,\al}$ assuming that $F$ is of class $(K)$, homogeneous of degree $1$, and concave. In addition there should exist a strictly convex function $\psi\in C^2(\bar\Om)$. The ambient space was an arbitrary Riemannian manifold $N$, $\bar\Om$ was supposed to be compact, and should be covered by a normal Gaussian coordinate system $(x^\al)$.

All hypotheses are satisfied in the present situation: $\bar \Om$ is compact, the normal Gaussian coordinate system is given by choosing geodesic polar coordinates with center in $x_0$, the strictly convex function $\psi$ can be defined by
\begin{equation}
\psi=\tfrac12\abs{x^0}^2,
\end{equation}
where $x^0$ is the radial distance to $x_0$, as one easily checks observing that the level hypersurfaces $\{x^0=\const\}$ which intersect $\bar\Om$ are all uniformly strictly convex, and $F$ is homogeneous of degree $1$ and therefore also concave.

\nocite{nk,evans:krylov,schneider:book}

\bibliographystyle{hamsplain}
\providecommand{\bysame}{\leavevmode\hbox to3em{\hrulefill}\thinspace}
\providecommand{\href}[2]{#2}



\end{document}